\documentclass[12pt]{article}
\usepackage{amsmath,amssymb,amscd}
\newtheorem{theorem}{Theorem} 
\newtheorem{lemma}{Lemma}

\newcommand{\qed}{{\unskip\nobreak\hfil\penalty50\hskip2em\vadjust{}
\nobreak\hfil$\Box$\parfillskip=0pt\finalhyphendemerits=0\par}}

\renewcommand{\P}{(P^{-1}(0))}
\newcommand{\z}{z^{(1)}}
\newcommand{\m}{multi}
\renewcommand{\o}{polynomial}
\newcommand{\cd}{\cdot}
\newcommand{\bs}{\bigskip}
\newcommand{\A}{{\bar A}(\theta)}
\renewcommand{\th}{\theta}

\begin{document}
\setlength{\baselineskip}{18pt}

\begin{center}{\Large\bf Multivariate Gauss-Lucas Theorems}\\ \bigskip\bigskip
{\bf Marek Kanter}\end{center}

\bs\bs
\begin{quote} {\bf Abstract.} A multivariate Gauss-Lucas theorem is proved, sharpening
and generalizing previous results on this topic. The theorem is stated in terms of a seemingly
new notion of convexity. Applications to \m variate stable polynomials are given.
\end{quote}

\bs\bs\bs

Let $P(z)$ be a polynomial in the complex variable $z$ with non-null derivative $Q(z)$.
The classical Gauss-Lucas theorem is the relation
\begin{equation}\label{1}
Q^{-1}(0)\subset H\P
\end{equation}
where $H\P$ is the convex hull of the roots of $P(z)$ in the complex plane.
(See Marden [5, Theorem 6.1].) It is intriguing to look for a \m dimensional extension
of this elegant result.

To establish notation, let $P(z)$ stand for a \o \ in the $M$ complex variables
$z\!=\!(z_1,\dots ,z_M)$. Given $P(z)$, let $Q(z)$ stand for the partial derivative of $P(z)$
with respect to $z_k$. The desired extension of (\ref{1}) has form
\begin{equation}\label{2}
Q^{-1}_k(0)\subset H\P \ ,
\end{equation}

\vfill\noindent
{\footnotesize {\bf 2000 Mathematics Subject Classification.}\\
Primary \ 12D10, \ 26C10, \ 30C15\\
Secondary \ 32A60, \ 26B25
Gauss-Lucas theorem, convex hull, separately convex, stable polynomials.}

\newpage\noindent
under the condition that $Q_k(z)$ is non-null. The convex hull in (\ref{2}) is taken with
respect to the usual convex structure on $C^M$ treated as isomorphic to $R^{2M}$. ($C$ stands for
the complex numbers and $R$ stands for the real numbers.)

It is somewhat curious that no proof of (\ref{2}) seems to be available in the literature
under the conditions stated. The papers [3] and [4] prove a direct generalization of
(\ref{1}) in the form
\begin{equation}\label{3}
\bigcap^M_{k=1} Q^{-1}_k(0)\subset H\P \ ,
\end{equation}
under the condition that $P(z)$ has finitely many zeros. (The left-hand side of (\ref{3})
represents the critical points of $P(z)$ in $C^M$.)

The finiteness restriction on the null set of $P(z)$ in (\ref{3}) is removed
in a related result about multivariate stable \o s, i.e.~\o s that have no roots
$z\!=\!(z_1,\dots ,z_M)$ with $Im(z_k)>0$ for all $k$. It is shown in 
[1, Theorem 3.1] and [7, Lemma 2.4] that the $k$-th partial derivative $Q_k(z)$ of
such a \o \ $P(z)$ is also stable, a result called the ``multivariate Gauss-Lucas
theorem", even though the roots of $Q_k(z)$ are not shown to lie in $H\P$. This further refinement
is provided in Theorem 1, using a weaker notion of convexity than that used to
define $H\P$, and thereby obtaining a smaller convex hull for $P^{-1}(0)$.

Preparatory to the proof of Theorem 1, it is necessary to establish some
terminology. For $z\in C^M$ with $z\!=\!(z_1,\dots ,z_M)$, the variables $z_k$ are
called the coordinates of $z$. The notation $z^{(k)}$ stands for the $M\!-\!1$ coordinates
of $z$ obtained by omitting the coordinate $z_k$. Given a subset $A$ of $C^M$, the section
of $A$ determined by $z^{(k)}$ is defined by
\[
A(z^k)\equiv \{ y_k\in C: y\in A, \ y^{(k)}\!=\!z^{(k)}\} \ .
\]

A set $A\subset C^M$ is called {\it separately convex in} $C^M$ if $A(z^{(k)})$ is a
convex subset of the complex plane for all $k\in \{1,\dots ,M\}$ and all $z\in C^M$.
If two sets are separately convex in $C^M$, then so is their intersection. This allows
the definition of the convex hull of a set $A\subset C^M$ as the smallest set
which is separately convex in $C^M$ and contains $A$, denoted $H_2(A)$.

\begin{theorem} If $P(z)$ is a multivariate \o \ defined on $C^M$,
with non-null $k$-th partial derivative $Q_k(z)$, then 
\begin{equation}\label{4}
 Q^{-1}_k(0)\subset H_2\P \ .
\end{equation}
\end{theorem}

{\bf Proof.} Without loss of generality let $k\!=\!1$. Suppose $z\in C^M$ with
$Q_1(z)=0$. To prove the theorem it suffices to show that
$z\in H_2\P$.

Let $f(w)=P(w,z^{(1)})$ for $w\in C$. Regard $f(w)$ as a \o \ in $w$ whose
coefficients depend on the fixed vector $z^{(1)}$. Note that
\[
f'(z_1)=Q_1(z_1,\z )=Q_1(z) =0 \ ,
\]
by assumption. The univariate Gauss-Lucas theorem implies that
$z_1\in H(f^{-1}(0))$, where $H(f^{-1}(0))$ stands for the usual convex hull
in $C$ of the roots of $f(w)$. Let $A$ be any subset of $C^M$ which is separately
convex in $C^M$ and contains $P^{-1}(0)$. Then the section $A(\z)$ is a convex subset of $C$
containing $f^{-1}(0)$; this implies that $A(\z)$ contains $H(f^{-1}(0))$.
In particular, $z_1\in A(\z)$, hence $z\in A$. By choice of $A$, this shows
that $z\in H_2\P$ and proves the theorem. \qed

\bs
The application of Theorem 1 to
multivariate stable \o s can easily be broadened by generalizing the definition of stability.
Let $\th=(\th_1,\dots ,\th_M)$ be a real vector and define
\[
A(\th) = \{z\in C^M: Im(e^{i\th_k}\,z_k)>0 \ {\mbox{ for }} 1\leq k\leq M\} \ .
\]
Let $\A$ stand for the complement of $A(\th)$ in $C^M$.

\begin{lemma} $\A$ is separately convex in $C^M$.\end{lemma}

{\bf Proof.}  It suffices to prove that for $z\in C^M$ the section of $\A$ determined
by $z^{(k)}$ is a convex subset of the complex plane for all $k$. Without loss of generality
set $k\!=\!1$ and let $z_j=x_j+iy_j$, where $x_j$ and $y_j$ are real. Let $c$ stand
for the minimum of $Im(e^{i\th_j}\,z_j)$ for $j\geq 2$. The section of $\A$ 
determined by $z^{(1)}$ is the set 
\[
\{x+iy\in C: (x\sin\th_1 + y\cos\th_1) \wedge c\leq 0\} \ .
\]
It is easy to see that this set is convex by considering the cases $c>0$ and $c\leq 0$
separately. \qed

\begin{theorem} Suppose the multivariate \o \ $P(z)$ is $\th$-stable, meaning
that it has no zero in $A(\th)$. Then any  non-null partial derivative
$Q_k(z)$ is also $\th$-stable.
\end{theorem}

{\bf Proof.} Apply Lemma 1 to get $H_2\P\subset\A$ and finish using Theorem 1. \qed

\bs It is clear that $H_2\P \subset H\P$, hence (\ref{4}) sharpens (\ref{2}). With
a view towards further sharpening (\ref{2}), it is natural to compare separate
convexity in $C^M$ with the analogous notion of separate convexity in $R^N$.
The latter notion can be defined by substituting $R$ for $C$ in the definition of separate
convexity in $C^M$ and has been studied under many names. (See [2] and [6].)

If $C^M$ is identified with $R^{2M}$ via the real and imaginary parts of the
complex coordinates in $C^M$, then any separately convex subset of $C^M$ is
separately convex in $R^{2M}$. Letting $H_1\P$ stand for the smallest separately
convex subset of $R^{2M}$ containing the roots of $P(z)$, it follows that
\begin{equation}\label{5}
H_1\P \subset H_2\P \ .
\end{equation}
However, (\ref{4}) is not valid if $H_1\P$ is substituted for $H_2\P$, even if $M\!=\!1$
(Example 1). Thus, no straightforward sharpening of (\ref{4}) seems to be available.

\bs{\bf Example 1.} When  $M\!=\!1$, the univariate case, conditions are given
for cubic and quadratic \o s $P(z)$ so that (\ref{4}) is not valid if $H_2\P$ is
replaced with $H_1\P$.

If $P(z)$ is a cubic \o \ with real coefficients and real roots, then $H_1\P$ is
equal to $H_2\P$, so the substitution does not affect the validity of (\ref{4}).
Therefore, assume that $P(z)$ has roots $a\!+\!bi$, $a\!-\!bi$, and $c$
(with $a,b,c$ real and $b\neq 0$). Setting the coefficient of $z^3$ in $P(z)$ to be 1,
yields that
\[
P(z) = (z-c)((z-a)^2+b^2) \ .
\]
The derivative $Q(z)=P'(z)$ is calculated to be
\[
Q(z) = 3z^2-(4a+2c)z+ a^2+b^2+2a c \ .
\]
The set $H_1\P$ is connected, consisting of the line segment joining
$a+bi$ to $a-bi$ and the line segment joining $a$ to $c$.

The roots of $Q(z)$ can be written as 
\[
w_i = \frac{2a}{3} + \frac c3 \pm \frac i3\, (3b^2-(a-c)^2)^{\frac 12} \ ,
\]
where $3b^2 >(a-c)^2$ by the assumption that not all roots are real. If $a\neq c$,
then it is clear that the roots of $Q(z)$ do not lie in the set $H_1\P$.
If $a\!=\!c$, then $H_1\P$ consists of the line segment joining $a+bi$ to $a\!-\!bi$
and the roots of $Q(z)$ are of the form
$a\pm \big( \frac{i}{\sqrt{3}}\big)b\in H_1\P$.

In summary, if $P(z)$ is a cubic \o \ with real coefficients, then the
roots of $Q(z)$ lie in $H_1\P$ if and only if the roots of $P(z)$ lie in a
straight line parallel to either the real or imaginary axis in the complex plane.
If $P(z)$ is a quadratic \o \ with complex coefficients, then the same result
is true. (Note that in the quadratic case, $H_1\P$ consists of the two
roots of $P(z)$ if these roots do not lie on a line parallel to either the real
or imaginary axis. Thus $H_1\P$ is not connected in this case.)


\bs\begin{center}{\bf References}\end{center}

\begin{description}
\item{[1]} J.~Borcea, P.~Branden: {\it The Lee-Yang and Polya-Schur programs. II.
Theory of stable \o s and applications}, Comm. Pure Appl. Math.  62 (2009),
1595--1631.

\item{[2]} B.~Dacorogna: {\it Direct methods in the calculus of variations.}
Second edition, Applied Mathematical Sciences, vol.~78, Springer, New York, 2008.

\item{[3]} J.~B.~Diaz, D.~B.~Shaffer: {\it A generalization to higher
dimensions of a theorem of Lucas concerning the zeros of the derivative of
a \o \ of one complex variable}, Appl. Anal.  6 (1977), 109--117.

\item{[4]} A.~W.~Goodman: {\it Remarks on the Gauss-Lucas theorem in higher
dimensional space}, Proc. Amer. Math. Soc. 55 (1976), 97--12.

\item{[5]} M.~Marden: {\it Geometry of polynomials}, Math. Surveys 3.
Second edition, Amer. Math. Soc., Providence, RI (1966).

\item{[6]} T.~Ottman, E.~Soisulon-Soisinen, D.~Wood: {\it On the definition and
computation of rectilinear convex hulls}, Information Sci.  33 (1984),
157--171.

\item{[7]} D.~G.~Wagner: {\it Multivariate stable \o s: Theory and Applications},
Bull. Amer. Math. Soc.  48 (2011), 53--83.

\end{description}

\bs\bs\hfill{}\begin{tabular}{r}
{\sc 1216 Monterey Avenue, Berkeley, California 94707}\\
{\it E-mail address}  \ {\tt mrk$@$cpuc.ca.gov}
\end{tabular}

\end{document}